\documentclass[11pt]{amsart}
\usepackage{amssymb}
\usepackage{amsmath}
\usepackage{amsthm}
\usepackage{latexsym}
\usepackage{amsfonts}
\usepackage{graphicx}
\usepackage{graphics}
\usepackage[T1]{pbsi}

\newtheorem{thm}{Theorem}[section]   
\newtheorem{lem}{Lemma}
\newtheorem{cor}{Corollary}
\newtheorem{Def}{Definition}

\theoremstyle{definition}
\newtheorem*{Proof}{Proof}

\newcommand{\dis}{\displaystyle}
\newcommand{\ra}{\;\rightarrow\;}

\newcommand{\bi}{\beta}
\newcommand{\ga}{\gamma }
\newcommand{\Ga} {{\varGamma}}
\newcommand{\Sig} {{\varSigma}}

\newcommand{\De} {{\varDelta}}

\newcommand{\f}{\varphi}

\newcommand{\la}{\lambda }
\newcommand{\mi}{\mu }

\newcommand{\si}{\sigma }

\newcommand{\ti}{\tau }

\newcommand{\oo}{\omega}

\newcommand{\R}{\mathbb{R}}
\newcommand{\Z}{\mathbb{Z}}
\newcommand{\N}{\mathbb{N}}

\newcommand{\ssum}{\sum\limits}

\newcommand{\cf}{{\mathcal{F}}}
\newcommand{\ca}{{\mathcal{A}}}

\newcommand{\ct}{{\mathcal{T}}}

\newcommand{\cm}{{\mathcal{M}}}

\newcommand{\ld}{\ldots}

 \newcommand{\loc}{\mbox{\footnotesize loc}}
\newcommand{\hs}{\hfill$\square$}

\begin{document}

\title[DYADIC MAXIMAL OPERATORS]{Extremal sequences for the Bellman function of the dyadic maximal operator}
\author{Eleftherios N. Nikolidakis}

\footnotetext{\hspace{-0.5cm}2010 MCS: 42B25} \footnotetext{\hspace{-0.5cm}Keywords: Bellman, dyadic, maximal}
\footnotetext{\hspace{-0.5cm}This research has been cofinanced by the European and Greek national funds through the operational program
Education and Lifelong Learning" of the National Strategic Reference Framework (NSRF), Aristeia code: MAXBELLMAN 2760, Research code: 70/3/11913}
\footnotetext{\hspace{-0.5cm} E-mail address: lefteris@math.uoc.gr}
\date{}
\maketitle

\noindent
{\bf Abstract:} We give a characterization of the extremal sequences for the Bellman function of the dyadic maximal operator.\ In fact we prove that they behave approximately like eigenfunctions of this operator for a specific eigenvalue.
\section{Introduction}\label{sec1}
\noindent

The dyadic maximal operator on $\R^n$ is a usefull tool in analysis and is defined by
\begin{eqnarray}
\cm_d\phi(x)=\sup\bigg\{\frac{1}{|Q|}\int_Q|\phi(u)|du:x\in Q,\;Q\subseteq\R^n\;\text{is a dyadic cube}\bigg\}  \label{eq1.1}
\end{eqnarray}
for every $\phi\in L^1_{\loc}(\R^n)$, where $|\cdot|$ is the Lesbesgue measure on $\R^n$ and the dyadic cubes are those formed by the grids $2^{-N}\Z^n$, $N=0,1,2,\ld\,.$

It is well known that it satisfies the following weak type (1,1) inequality
\begin{eqnarray}
|\{x\in\R^n:\cm_{d}\phi(x)>\la\}|\le\frac{1}{\la}\int_{\{\cm_d
\phi>\la\}}|\phi(u)|du, \label{eq1.2}\mathcal{}
\end{eqnarray}
for every $\phi\in L^1(\R^n)$ and $\la>0$.

From (\ref{eq1.2}) it is not difficult to prove the following $L^p$-inequality
\begin{eqnarray}
\|\cm_d\phi\|_p\le\frac{p}{p-1}\|\phi\|_p, \ \ \label{eq1.3}
\end{eqnarray}
for every $p>1$ and $\phi\in L^p(\R^n)$,
and this can be done by using the well known Doob's method for the dyadic maximal operator.

It is also easy to see that (\ref{eq1.2}) is best possible, while (\ref{eq1.3}) is also best possible as can be seen in \cite{Wang} (see \cite{Burka} and \cite{Burkb} for general martingales).

Our aim in this article is to study this maximal operator and one way to do this is to find certain refinements of the inequalities satisfied by it such as (\ref{eq1.2}) and (\ref{eq1.3}). Concerning (\ref{eq1.2}) refinements have been made in \cite{Nikol1}, \cite{Nikol3} and \cite{Nikol2}. Refinements of (\ref{eq1.3}) can be found in \cite{Mel1} or even more general in \cite{Mel2}.

In order to refine (\ref{eq1.3}) we should introduce the following function
\begin{eqnarray}
B^Q_p(f,F)=\sup\bigg\{\frac{1}{|Q|}\int_Q(\cm_d\phi)^p:\phi\ge0,\; Av_{Q}(\phi)=f, \; Av_{Q}(\phi^p)=F\bigg\}  \label{eq1.4}
\end{eqnarray}
where, $p>1$, $0<f^p\le F$, $Q$ is a fixed dyadic cube in $\R^n$, $\phi\in L^p(Q)$ and
\[
Av_{Q}(h)=\frac{1}{|Q|}\int_Q|h(u)|du,
\]
for every $h\in L^1(Q)$.
This is the so-called Bellman function of two variables associated to the dyadic maximal operator. Thus by considering the above function we refine
(\ref{eq1.3}), by adding a norm variable, which is the $L^1$-norm of $\phi$, and which we consider to be equal to a fixed constant $f$.

In fact this function has been explicitly computed. Actually this is done in a much more general setting of a non-atomic probability measure space $(X,\mi)$, where the dyadic sets are now given in a family of sets $\ct$ (called tree), which satisfies conditions similar to those that are satisfied by the dyadic cubes on $[0,1]^n$ (for details see section 2). We then define the associated dyadic maximal operator $\cm_\ct$ by
\begin{eqnarray}
\cm_\ct\phi(x)=\sup\bigg\{\frac{1}{\mi(I)}\int_I|\phi|d\mi:x\in I\in\ct\bigg\},  \label{eq1.5}
\end{eqnarray}
for every $\phi\in L^1(X,\mi)$.

 The Bellman function of two variables for $p>1$ associated to $\cm_\ct$ is now given by
\begin{eqnarray}
\hspace*{1,5cm} B^{\ct}_p(f,F)=\sup\bigg\{\int_X(\cm_\ct\phi)^pd\mi:\phi\ge0,\;\int_X\phi d\mi=f,\;\int_X\phi^pd\mi=F\bigg\},  \label{eq1.6}
\end{eqnarray}
where $0<f^p\le F$.

In \cite{Mel1}, (\ref{eq1.6}) has been found equal to $F\oo_p(f^p/F)^p$ where $\oo_p:[0,1]\longrightarrow\Big[1,\dfrac{p}{p-1}\Big]$ is the inverse function $H^{-1}_p$ of $H_p$ defined for $z\in\Big[1,\dfrac{p}{p-1}\Big]$ by $H_p(z)=-(p-1)z^p+pz^{p-1}$.
This gives us as an immediate consequence that it is independent of the measure space $(X,\mi)$ and the tree structure of $\ct$.

For the evaluation of this function the author in \cite{Mel1} introduced a technique which enabled him to compute it. This is based on an effective linearization of the dyadic maximal operator that holds for an adequate set of functions, called $\ct$-good. Certain sharp inequalities were proved in \cite{Mel1} by using H\"{o}lder's inequality upon suitable subsets of $X$ in an effective way. After the evaluation of (\ref{eq1.6}) he was also able to evaluate other more general Bellman functions of $\cm_\ct$ that involve three parameters. The evaluations of these new Bellman functions, which are connected with the Dyadic Carleson Imbedding Theorem and others, are based on the result of (\ref{eq1.6}) entirely and are proved by its application on certain elements of the tree $\ct$.

The next step for studying the dyadic maximal operator is to investigate the opposite problem for the Bellman function related to Kolmogorov's inequality which has been worked out in \cite{Mel3}. More precisely the following function
\begin{eqnarray}
\hspace*{1,5cm} B_q(f,h)=\sup\bigg\{\int_X(\cm_\ct\phi)^qd\mi:\phi\ge0,\;\int_X\phi d\mi=f,\;\int_X\phi^qd\mi=h\bigg\},  \label{eq1.7}
\end{eqnarray}
has been computed there, where $0<h\le f^q$ and $q\in (0,1)$ is a fixed constant.

In \cite{Mel3} the authors precisely computed the above function by using the linearization technique introduced in \cite{Mel1}. The situation is now different and new methods were found in order that (\ref{eq1.7}) be evaluated.

Now the following has been proved in \cite{Nikol7}. \vspace*{0.2cm}\\
\noindent
{\bf Proposition:}{\em Let $(\phi_n)_n$ be a sequence of nonnegative functions in $L^1(X,\mi)$ such that $\int_X\phi_nd\mi=f$ and
$\int_X\phi^p_nd\mi=F$ for all $n\in N$. If $(\phi_n)_n$ is extremal for (\ref{eq1.6}), then for every $I\in\ct$ we have that
 $\dis\lim_n \frac {1}{\mi(I)}\int_I\phi_nd\mi=f$ and $\dis\lim_n \frac {1}{\mi(I)}\int_I\phi_n^pd\mi=F$. Moreover
 $\dis\lim_n \frac {1}{\mi(I)}\int_I(\cm_\ct\phi_n)^pd\mi=B_p^{\ct}(f,F)$.} \vspace*{0.2cm}

This gives as an immediate result that there do not exist extremal functions for (\ref{eq1.7}). This is true because if $\ct$ differentiates $L^1(X,\mi)$ we would have for any extremal $\phi$ that it should be constant almost everywhere on $X$, so that $F=f^p$ which is a trivial case that
we do not consider.

Thus our interest is for those sequences of functions $(\phi_n)_n$ that are extremal for this Bellman function. That is $\phi_n:(X,\mi)\ra\R^+$, $n\in N$ must satisfy
\[
\int_X\phi_nd\mi=f,\;\int_X\phi^p_nd\mi=F \ \ \text{and} \ \ \lim_n\int_X(\cm_\ct\phi_n)^pd\mi=F\oo_p(f^p/F)^p.
\]
Our aim in this paper is to give a characterization of these extremal sequences of functions. For this reason we restrict ourselves to the class of  $\ct$-good functions, that is enough to describe the problem as it was settled in \cite{Mel1} (see section 3). We give now the statement of our main result.

\noindent
{\bf Theorem A:} {\em Let $(\phi_n)_n$ be a sequence of nonnegative, $\ct$-good functions such that $\int_X\phi_nd\mi=f$ and $\int_X\phi^p_nd\mi=F$. Then $(\phi_n)_n$ is extremal for (\ref{eq1.6}), if and only if $$\dis\lim_n\int_X|\cm_\ct\phi_n-c\phi_n|^pd\mi=0,$$ for $c=\oo_p(f^p/F)$.} \vspace*{0.2cm}

That is $(\phi_n)_n$ is an extremal sequence for (\ref{eq1.6}), if and only if its terms behave approximately, in $L^p$, like eigenfunctions of $\cm_\ct$,
for the eigenvalue $c=\oo_p(f^p/F)$.

For the proof of the above theorem we use the technique introduced in \cite{Mel1} for the evaluation of (\ref{eq1.6}), which we generalize
in two directions (see theorems 3.1 and 3.2) and by using these we prove theorem 3.3 for the extremal sequences we are interested in. This theorem is in fact a weak form of theorem A. It is proved by producing two inequalities that involve the $L^p$-integrals of $\cm_\ct\phi$ and $\phi$ over measurable subsets $A\subset X$ that have a certain form with respect to the tree $\ct$ and the function $\phi$. More precisely A is a union of certain elements of $S_\phi$ or a complement of such a set, where $S_\phi$ is a subtree of $\ct$ that depends on $X$ and gives all the information we need for $\cm_\ct\phi$ (for the definition of $S_\phi$ see section 2). Using theorems 3.1 and 3.2 we eventually reach to theorem 3.3.

In order to prove theorem A we need to apply theorem 3.3 to a new extremal sequence $(g_{\phi_n})$ which is arbitrarily close to $(\phi_n)_n$ in the $L^p$ sense. In fact $g_{\phi_n}$ is defined properly  on suitable subsets of $X$ where $\phi_n$
is defined. The number of different values of $g_{\phi_n}$ on each of these subsets are at most two with the one being zero. Then we prove that
the measure of the set where $g_{\phi_n}$ is zero tends to zero by using the fact that $(g_{\phi_n})$ is extremal sequence for (\ref{eq1.6}). Thus we can arrange everything so that this new extremal sequence is constant on those suitable sets. We denote this new sequence by $(g_{\phi_n}')$. Since $g_{\phi_n}'$ is constant on each one of the suitable subsets of $X$, we are in position to apply theorem 3.3 to it and by using some additional technical lemmas we finally reach to theorem A.

We should also note that additional work concerning the Bellman functions and certain symmetrization principles for the dyadic maximal operator can be seen in \cite{Mel2} and \cite{Nikol4}. It is also worth saying that in \cite{Nikol6} it has been given an alternative method for the evaluation of the Bellman function (\ref{eq1.6}). Also we need to remind that the phenomenon that the norm of a maximal operator is attained by a sequence of eigenfunctions of such a maximal operator can be seen in \cite{Graf} and \cite{Colz}. So by considering the results of this paper one might guess that it shouldn't be rare and and may occur in other settings also, such as square functions or other dyadic operators. Finally we mention that the extremizers for the Bellman function of three variables related to Kolmogorov's inequality have been characterized in \cite{Nikol5}.

\section{Preliminaries}\label{sec2}
\noindent

Let $(X,\mi)$ be a non-atomic probability measure space. We give the following from \cite{Mel1}.
\begin{Def}\label{Def2.1}
 A set $\ct$ of measurable subsets of $X$ will be called a tree if the following are satisfied
\begin{enumerate}
\item[i)] $X\in\ct$ and for every $I\in\ct$, $\mi(I)>0$.
\item[ii)] For every $I\in\ct$ there corresponds a finite or countable subset $C(I)$ of $\ct$ containing at least two elements such that
\begin{itemize}
\item[a)] the elements of $C(I)$ are pairwise disjoint subsets of $I$
\item[b)] $I=\cup C(I)$.
\end{itemize}
\item[iii)] $\ct=\bigcup\limits_{m\ge0}\ct_{(m
)}$, where $\ct_{(0)}=\{X\}$ and $\ct_{(m+1)}=\bigcup_{I\in\ct_{(m)}}C(I)$.
\item[iv)] The following holds
\[
\lim_{m\ra\infty}\sup_{I\in\ct_{(m)}}\mi(I)=0.
\]

\end{enumerate}
\end{Def}
The following is presented in \cite{Mel1}, and is a consequence of the properties i)-iv) of Definition 2.1, that a tree $\ct$ satisfies.
\begin{lem}\label{lem2.1}
For every $I\in\ct$ and every $a\in(0,1)$ there exists a subfamily $\cf(I)\subseteq\ct$ consisting of pairwise disjoint subsets of $I$ such that
\[
\mi\bigg(\bigcup_{J\in\cf(I)}J\bigg)=\sum_{J\in\cf(I)}\mi(J)=(1-a)\mi(I).
\]
\end{lem}

Now given a tree $\ct$ we define the maximal operator associated to it as follows
\[
\cm_\ct\phi(x)=\sup\bigg\{\frac{1}{\mi(I)}\int_I|\phi|d\mi:\;x\in I\in\ct\bigg\},  \ \
\]
for every $\phi\in L^{1}(X,\mi)$.
Then one can see in \cite{Mel1}, the following.
\begin{thm}\label{thm2.1}
The following equality is true
\[
\sup\bigg\{(\cm_\ct\phi)^pd\mi:\;\phi\ge0,\;\int_X\phi d\mi=f,\;\int_X\phi^pd\mi=F\bigg\}=F\oo_p(f^p/F)^p,
\]
for every $f,F$ such that $0<f^p\le F$.
\end{thm}

Additionally we give the notion of the extremal sequence as
\begin{Def}\label{Def2.2}
Let $(\phi_n)_n$ be a sequence of $\mi$-measurable nonnegative functions defined on $X$, $p>1$ and $0<f^p\le F$. Then $(\phi_n)_n$ is called $(p,f,F)$ extremal or simply extremal if the following hold:
\[
\int_X\phi_nd\mi=f,\;\int_X\phi^p_nd\mi=F \ \ \text{and} \ \ \lim_n\int_X(\cm_\ct\phi_n)^pd\mi=F\oo_p(f^p/F)^p.
\]
\end{Def}
\section{Characterization of the extremal sequences}\label{sec3}
\noindent

We describe now the effective linearization for the operator $\cm_\ct$ that was introduced in \cite{Mel1} which is valid for certain class of functions $\phi$.

For every $\phi\in L^1(X,\mi)$ nonnegative and $I\in\ct$ we define $Av_{I}(\phi)=\dfrac{1}{\mi(I)}\int_I\phi d\mi$.

We will say that $\phi$ is $\ct$-good if the set
\[
\ca_\phi=\{x\in X:\cm_\ct\phi(x)>Av_{I}(\phi)\ \ \text{for all} \ \ I\in\ct \ \ \text{such that} \ \ x\in I\}
\]
has $\mi$-measure zero.

Let now $\phi$ be $\ct$-good and $x\in X\setminus\ca_\phi$.
We define $I_\phi(x)$ to be the largest in the nonempty set
\[
\{I\in\ct:\;x\in I \ \text{and} \ \ \cm_\ct\phi(x)=Av_{I}(\phi)\}.
\]
Suppose now that $I\in\ct$. We define the following
\[
A(\phi,I)=\{x\in X\setminus\ca_\phi:\;I_\phi(x)=I\}\subseteq I,
\]
\[
S_\phi=\{I\in\ct:\;\mi(A(\phi,I))>0\}\cup\{X\}.
\]
Obviously then $\cm_\ct\phi=\sum\limits_{I\in S_\phi}Av_{I}(\phi)\chi_{A(\phi,I)}$, $\mi$-a.e., where $\chi_E$ is the characteristic function of $E$.

We define also the following correspondence $I\ra I^\ast$ by: $I^\ast$ is the smallest element of $\{J\in S_\phi:\;I\subsetneq J\}$. It is defined for every $I\in S_\phi$ except $X$. Then it is obvious that the $A(\phi,I)$ are pairwise disjoint and that $\mi\Big(\bigcup\limits_{I\notin S_\phi}(A(\phi,I))\Big)=0$, so that $\bigcup\limits_{I\in S_\phi}A(\phi,I)\approx X$, where by $A\approx B$ we mean that $\mi(A\setminus B)=\mi(B\setminus A)=0$.

Now the following is a consequence of the above.
\setcounter{lem}{0}
\begin{lem}\label{lem3.1}
Let $\phi$ be $\ct$-good and let also $I\in \ct$, $I\neq X$. Then $I\in S_\phi$ if and only if for every $J\in \ct$ that contains properly $I$ we have that
$Av_J(\phi)<Av_I(\phi)$.
\end{lem}
\begin{Proof}
Suppose that $I\in S_\phi$. Then $\mi(A(\phi,I))>0$. Thus $A(\phi,I)\neq\emptyset$, so there exists $x\in A(\phi,I)$. By the definition of
$A(\phi,I)$ we have that $I_{\phi}(x)=I$, that is $I$ is the largest element of $\ct$ such that $\cm_\ct\phi(x)=Av_I(\phi)$. As a consequence the implication
stated in our lemma holds.

Conversely suppose that $I\in \ct$ and for every $J\in \ct$ that contains properly $I$ we have that $Av_J(\phi)<Av_I(\phi)$.
Then since $\phi$ is $\ct-good$, we have that for every $x\in I\setminus\ca_\phi$ there exists $J_x=I_{\phi}(x)$ in $S_\phi$ such that  $\cm_\ct\phi(x)=Av_{J_x}(\phi)$ and
$x\in J_x$. By our hypothesis we must have that $J_x\subseteq I$. Consider the family $S^1=(J_x)_{x\in I\setminus\ca_\phi}$. This obviously has the following property: $\bigcup\limits_{x\in I\setminus\ca_\phi}J_x\approx I$. Choose now a pairwise disjoint subfamily $S^2=(J_i)_i$ with $X\approx \cup J_i$. For this choice we just need to consider those $J_x\in S^1$ maximal under $\subseteq$ relation. Then by our construction
$Av_{J_i}(\phi)\geq Av_{I}(\phi)$. Suppose now that $I \notin S_\phi$. This means that $\mi(A(\phi,I))=0$, that is we must have for every $x\in I\setminus\ca_\phi$ that $J_x \subsetneq I$. Since $J_x$ belongs to $S_\phi$ for every such $x$, by the first part of the proof of this Lemma we conclude that
$Av_{J_x}(\phi)> Av_{I}(\phi)$ and as a consequence we have that $Av_{J_i}(\phi)> Av_{I}(\phi)$ for every $i$. Since $S^2$ is a decomposition of $I$, and because of the last mentioned inequality we reach to a contradiction. In this way we derive the proof of our lemma.  \hs
\end{Proof}
Now the following is proved in \cite{Mel1}.
\begin{lem}\label{lem3.2}
Let $\phi$ be $\ct$-good
\begin{enumerate}
\item[i)] If $I$, $J\in S_\phi$ then either $A(\phi,J)\cap I=\emptyset$ or $J\subseteq I$.
\item[ii)] If $I\in S_\phi$ then there exists $J\in C(I)$ such that $J\notin S_\phi$.
\item[iii)] For every $I\in S_\phi$ we have that
\[
I\approx\bigcup_{J\in S_\phi\atop J\subseteq I}A(\phi,J).
\]
\item[iv)] For every $I\in S_\phi$ we have that
\[
A(\phi,I)=I\setminus\bigcup_{J\in S_\phi\atop J^\ast\in I}J, \ \ \text{so that}
\]
\[
\mi(A(\phi,I))=\mi(I)-\sum_{J\in S_\phi\atop J^\ast=I}\mi(J). \vspace*{-0.5cm}
\]
\end{enumerate}
\end{lem}

From all the above we see that
\[
Av_{I}(\phi)=\frac{1}{\mi(I)}\sum_{J\in S_\phi\atop J\subseteq I} \;\int_{A(\phi,J)}\phi d\mi=:y_I
\]
where $I\in S_\phi$, and for those $I$ we also define
\[
x_I=a_I^{-1+\frac{1}{p}}\int_{A(\phi,I)}\phi d\mi,  \ \ \text{where} \ \ a_I=\mi(A(\phi,I)).
\]
We prove now the following
\begin{thm}\label{thm3.1}
Let $\phi$ be $\ct$-good function such that $\int\limits_
X\phi d\mi=f$. Let also $B=\{I_j\}$ be a family of pairwise disjoint elements of $S_\phi$, which is maximal on $S_\phi$ under $\subseteq$ relation.\ That is if $I\in S_\phi$ then $I\cap(\cup I_j)\neq\emptyset$. Then the following inequality holds
\[
\int_{X\setminus\bigcup\limits_jI_j}\phi^pd\mi\ge\frac{f^p-\ssum_j\mi(I_j)y^p_{I_j}}
{(\bi+1)^{p-1}}+\frac{(p-1)\bi}{(\bi+1)^p}\int_{X\setminus\bigcup\limits_jI_j}
(\cm_\ct\phi)^pd\mi
\]
for every $\bi>0$, where $y_{I_j}=Av_{{I_j}}(\phi)$.
\end{thm}
\begin{Proof}
We follow \cite{Mel1}. We obviously have that
\setcounter{equation}{0}
\begin{eqnarray}
\int_{X\setminus\cup I_j}\phi^pd\mi=\sum_{I\supsetneq\text{piece}(B)\atop I\in S_\phi}\int_{A(\phi,I)}\phi^pd\mi,  \label{eq3.1}
\end{eqnarray}
where by writing $I\supsetneq\text{piece}(B)$ we mean that $I\supsetneq I_j$ for some $j$. In fact (\ref{eq3.1}) is true since $X\setminus\bigcup\limits_jI_j\approx\bigcup\limits_{J\in S_\phi\atop I\supsetneq\text{piece}(B)}A(\phi,I)$ in view of the maximality of $B$ and Lemma \ref{lem3.2}.

Now from (\ref{eq3.1}) we have by H\"{o}lder's inequality that
\begin{eqnarray}
\int_{X\setminus\bigcup\limits_jI_j}\phi^pd\mi\ge\sum_{I\in S_\phi\atop I\supsetneq\text{piece}(B)}x^p_I=\sum_{I\in S_\phi\atop I\supsetneq\text{piece}(B)}\frac{\Big(\int\limits_{A(\phi,I)}\phi d\mi\Big)^p}{a_I^{p-1}}.  \label{eq3.2}
\end{eqnarray}
It is also true that
\[
\mi(I)y_I=\sum_{J\in S_{\phi}\atop J^\ast=I}\mi(J)y_J+\int_{A(\phi,I)}\phi d\mi, \ \ \text{for every} \ \ I\in S_\phi.
\]
Thus by using H\"{o}lder's inequality in the form
\[
\frac{(\la_1+\cdots+\la_m)^p}{(\si_1+\cdots+\si_m)^{p-1}}\le\frac{\la^p_1}{\si_1^{p-1}}+
\frac{\la_2^p}{\si^{p-1}_2}+\cdots+\frac{\la^p_n}{\si^{p-1}_m},
\]
we have
\begin{align}
\int_{X\setminus\cup I_j}\phi^pd\mi\ge\sum_{I\in S_\phi\atop I\supsetneq\text{piece{(B)}}}
\frac{\Big(\mi(I)y_I-\ssum_{J\in S_{\phi}\atop J^\ast=I}\mi(J)y_J\Big)^p}{\Big(\mi(I)-\ssum_{J\in S_{\phi}\atop J^\ast=I}
\mi(J)\Big)^{p-1}} \nonumber \\
\ge\sum_{I\in S_\phi\atop I\supsetneq\text{piece{(B)}}}\bigg\{\frac{(\mi(I)y_I)^p}{(\ti_I\mi(I))^{p-1}}-
\sum_{J\in S_{\phi}\atop J^\ast=I}\frac{(\mi(J)y_J)^p}{((\bi+1)\mi(J))^{p-1}}\bigg\},  \label{eq3.3}
\end{align}
where $\ti_I=(\bi+1)-\bi\rho{_I}$, $\rho_I=\dfrac{a_I}{\mi(I)}$, $\bi>0$.

Then by (\ref{eq3.3}) we have because of the maximality of B that
\begin{eqnarray}
\int_{X\setminus\bigcup\limits_jI_j}\phi^pd\mi\ge\sum_{I\in S_\phi\atop I\supsetneq\text{piece}(B)}\frac{\mi(I)y^p_I}{\ti_I^{p-1}}-\sum_{(\ast)}
\frac{\mi(I)y^p_I}{(\bi+1)^{p-1}},  \label{eq3.4}
\end{eqnarray}
where the summation in $(\ast)$ is extended to:

(a) $I\in S_\phi$: $I\supsetneq\text{piece}(B)$ with $I\neq X$  and (b) $I\in S_\phi$ is a piece of $B$ ($I=I_j$, for some $j$).

As a consequence we can write
\begin{align}
\int_{X\setminus\cup I_j}\phi^pd\mi\ge&\frac{y^p_x}{\ti^{p-1}_x}+\sum_{\footnotesize{\begin{array}{c}
                                                                 I\in S_\phi  \\ [-0.8ex]
                                                                  I\neq X \\ [-0.8ex]
                                                                  I\supsetneq\text{piece}(B) \end{array}}}
\frac{1}{\rho_I}\bigg(\frac{1}{\ti^{p-1}_I}\nonumber
-\frac{1}{(\bi+1)^{p-1}}\bigg)a_Iy^p_I-    \nonumber\\
&-\frac{1}{(\bi+1)^{p-1}}\sum_j\mi(I_j)y^p_{I_j}.  \label{eq3.5}
\end{align}
It is easy now to see that
\begin{eqnarray}
\frac{1}{(\bi+1-\bi x)^{p-1}}-\frac{1}{(\bi+1)^{p-1}}\ge\frac{(p-1)\bi x}{(\bi+1)^p},  \label{eq3.6}
\end{eqnarray}
for any $x\in [0,1]$, in view of the mean value theorem on derivatives.

Then by (\ref{eq3.5}) we immediately conclude that
\begin{align}
\int_{X\setminus\cup I_j}\phi^pd\mi\ge&\frac{y^p_X}{\ti^{p-1}_X}+\frac{(p-1)\bi}
{(\bi+1)^p}\sum_{\footnotesize{\begin{array}{c}
                                                                 I\neq X\\ [-0.8ex]
                                                                  I\in S_\phi \\ [-0.8ex]
                                                                  I\supsetneq\text{piece}(B) \end{array}}}
a_Iy^p_I-\frac{1}{(\bi+1)^{p-1}}\sum_j\mi(I_j)y^p_{I_j} \nonumber\\
=&\bigg[\frac{1}{((\bi+1)-\bi\rho_{X})^{p-1}}-\frac{(p-1)\bi\rho_{X}}
{(\bi+1)^p}\bigg]f^p+\frac{(p-1)\bi}{(\bi+1)^p}\sum_{I\in S_\phi\atop I\supsetneq\text{piece}(B)}a_Iy^p_I \nonumber \\
&-\frac{1}{(\bi+1)^{p-1}}\sum_j\mi(I_j)y^p_{I_j}.  \label{eq3.7}
\end{align}
On the other hand $\ssum_{I\in S_\phi\atop I\supsetneq\text{piece}(B)}a_Iy^p_I=\sum\limits_{X\setminus\cup I_j}(\cm_\ct\phi)^pd\mi$, so in view of (\ref{eq3.6}) we must have that
\[
\int_{X\setminus\cup I_j}\phi^p\ge\frac{f^p-\sum\mi(I_j)y^p_{I_j}}{(\bi+1)^{p-1}}+
\frac{(p-1)\bi}{(\bi+1)^p}\int_{X\setminus\cup I_j}(\cm_\ct\phi)^pd\mi,
\]
for every $\bi>0$, and the proof of the theorem is complete.  \hs
\end{Proof}

If we follow the same proof as above but now work inside any of the $I_j$ we obtain
\begin{thm}\label{thm3.2}
Let $\phi$ be $\ct$-good and $\ca=\{I_j\}$ be a pairwise disjoint family of elements of $S_\phi$. Then for every $\bi>0$ we have that:
\[
\int_{\bigcup\limits_jI_j}\phi^pd\mi\ge\frac{\sum\mi(I_j)y^p_{I_j}}{(\bi+1)^{p-1}}
+\frac{(p-1)\bi}{(\bi+1)^p}\int_{\bigcup\limits_jI_j}(\cm_\ct\phi)^pd\mi. \vspace*{-0.5cm}
\]
\end{thm}
Let us now prove the following generalization of theorem 3.1.
\begin{cor}\label{cor3.1}
Suppose that $\phi$ is $\ct$-good and $A=\{I_j\}$ be a pairwise disjoint family of elements of $S_\phi$. Then for every $\bi>0$
\[
\int_{X\setminus{\bigcup\limits_jI_j}}\phi^pd\mi\ge\frac{f^p-\ssum_j\mi(I_j)y^p_{I_j}}
{(\bi+1)^{p-1}}+\frac{(p-1)\bi}{(\bi+1)^p}\int_{X\setminus\bigcup\limits_jI_j}(\cm_\ct\phi)^pd\mi,
\]
where $f=\int\limits_X\phi d\mi$.
\end{cor}
\begin{Proof}
This is true since there exist families $B,\Ga$ of pairwise disjoint elements of $S_\phi$ with $B$ as in the statement of theorem \ref{thm3.1}, such that $B=\bigcup\limits_jI'_j$, $\Ga=\bigcup\limits_iJ_i$ with $\bigcup\limits_jI'_j=\Big(\bigcup\limits_jI_j\Big)\cup\Big(\bigcup\limits_iJ_i\Big)$ and the additional property that $I_j$ is disjoint to $J_i$ for every $j,i$.
Applying theorem \ref{thm3.1} for $B$ and theorem \ref{thm3.2} for $\Ga$ we obtain, by summing the respective inequalities, the proof of corollary \ref{cor3.1}.  \hs
\end{Proof}

As a consequence of the above we have the following.
\begin{thm}\label{thm3.3}
Let $(\phi_n)_n$ an extremal sequence consisting of $\ct$-good functions. Consider for every $n\in\N$ a pairwise disjoint family $\ca_n=\{I^n_j\}$ of elements of $S_{\phi_n}$ such that the following limit exists
\[
\lim_n\sum_{I\in\ca_n}\mi(I)y^p_{I,n}, \ \ \text{where} \ \ y_{I,n}=Av_{I}(\f_n), \ \ I\in\ca_n.
\]
Then
\[
\lim_n\int_{\cup\ca_n}(\cm\phi_n)^pd\mi=\oo_p(f^p/F)^p\lim_n\int_{\cup\ca_n}\phi^p_nd\mi
\]
meaning that if one of the limits on the above relation exists then the other also does and we have the stated equality.
\end{thm}
\begin{Proof}
In view of Theorem \ref{thm3.2} and Corollary \ref{cor3.1} we have that
\begin{eqnarray}
\hspace*{1cm}\int_{X\setminus\cup\ca_n}\phi^p_nd\mi\ge\frac{f^p-\ssum_{I\in\ca_n}\mi(I)y^p_{I,n}}{(\bi+1)^{p-1}}+
\frac{(p-1)\bi}{(\bi+1)^p}\int_{X\setminus\cup\ca_n}(\cm_\ct\phi_n)^pd\mi,   \ \ \text{and} \label{eq3.8}
\end{eqnarray}
\begin{eqnarray}
\hspace*{1cm}\int_{\cup\ca_n}\phi^p_nd\mi\ge\frac{\sum_{I\in\ca_n}\mi(I)y^p_{I,n}}{(\bi+1)^{p-1}}+
\frac{(p-1)\bi}{(\bi+1)^p}\int_{\cup\ca_n}(\cm_\ct\phi_n)^pd\mi,  \label{eq3.9}
\end{eqnarray}
for every $\bi>0$ and $n\in\N$.

Summing relations (\ref{eq3.8}) and (\ref{eq3.9}) for every $n\in\N$ we obtain
\begin{eqnarray}
F=\int_X\phi^p_nd\mi\ge\frac{f^p}{(\bi+1)^{p-1}}+\frac{(p-1)\bi}{(\bi+1)^p}\int_X
(\cm_\ct\phi_n)^pd\mi,  \label{eq3.10}
\end{eqnarray}
Since $(\phi_n)_n$ is extremal we have equality in the limit in (\ref{eq3.10}) for $\bi=\oo_p(f^p/F)-1$ (see \cite{Mel1}, relation (4.24)).

So we must have equality on (\ref{eq3.8}) and (\ref{eq3.9}) in the limit for this value of $\bi$.
Suppose now that $h_n=\ssum_{I\in\ca_n}\mi(I)y^p_{I,n}$ and that $h_n\ra h$. Now we can write (\ref{eq3.9}) in the form
\begin{eqnarray}
\int_{\cup_{\ca_n}}(\cm_\ct\phi_n)^pd\mi\le\bigg(1+\frac{1}{\bi}\bigg)
\frac{(\bi+1)^{p-1}\int\limits_{\cup\ca_n}\phi^p_nd\mi-h_n}{p-1},  \label{eq3.11}
\end{eqnarray}
(see \cite{Mel1}, relations (4.24) and (4.25)), for every $\bi>0$. The right hand side of (\ref{eq3.11}), $n\in\N$, is minimized for $\bi=\bi_n=\oo_p\Big(h_{n}\big/\int\limits_{\cup\ca_n}\phi^p_nd\mi\Big)-1$, as can be seen at the end of the proof of lemma 9 in \cite{Mel1},
or by making the related simple calculations.

Since, we have equality in the limit in (\ref{eq3.11}) we must have that
\begin{eqnarray}
\lim_n\frac{h_n}{\int\limits_{\cup\ca_n}\phi^p_nd\mi}=\frac{f^p}{F},  \label{eq3.12}
\end{eqnarray}
Thus (\ref{eq3.12}) and (\ref{eq3.11}) give
\[
\lim_n\int_{\cup\ca_n}(\cm_\ct\phi_n)^pd\mi=\oo_p(f^p/F)^p\lim_n\int_{\cup\ca_n}\phi^p_n
\]
and this holds in the sense stated above. This completes the proof of theorem 3.3. \hs
\end{Proof}

We need now some additional lemmas that we are going to state and prove below. First we prove the following.
\begin{lem}\label{lem3.3}
Let $\phi$ be $\ct$-good. Then we can associate to $\phi$, a measurable function defined on $X$, $g_{\phi}$, which attains two at most values
($c_J^{\phi}$ or 0) on certain subsets of $A(\phi,J)$, that decompose it, for every $J\in S_{\phi}$, and which is defined in a way that for every $I\in \ct$ which contains an element of $S_{\phi}$ (that is it is not contained in any of the $A_J$) we must have that $\int_Ig_{\phi}d\mi=\int_I\phi d\mi$.
Additionally for any $I\in S_{\phi}$ we have that $\int_{A_I}g^p_{\phi}d\mi=\int_{A_I}\phi^p d\mi$ and $\mi(\{\phi=0\}\cap A_I)\le \mi(\{g_{\phi}=0\}\cap A_I)$.
\end{lem}
\begin{Proof}
We define $g_{\phi}$ inductively using Lemma 3.2. Note that $A(\phi,X)=A_X= X\setminus\cup_{I\in S_{\phi}, I^{\ast}=X}I$.
We define first a function $g_{\phi}^{(1)}:X\rightarrow\R^+$ such that the integral relation mentioned above holds for this function and additionally
$g_{\phi}^{(1)}/A_X$ attains at most two values on certain subsets of $A_X$, which are in fact unions of elements of $\ct$, and which decompose
$A_X$. For this construction we proceed as follows. We set $g_{\phi}^{(1)}(x)=\phi(x)$, for $x\in X\setminus A_X$. We write
$A_X=\cup_{j}I_{j,X}$, where $(I_{j,X})_j$  is a family of elements of $\ct$, maximal with respect to the relation $I_{j,X}\subseteq A_X$.
For every $I_{j,X}$ there exists an integer $k_j>0$, such that $I_{j,X}\in \ct_{(k_j)}$. Then we consider the unique $I_{j,X}^{'}$ such that
$I_{j,X}\in C(I_{j,X}^{'})$, that is $I_{j,X}^{'}\in \ct_{(k_j-1)}$ and $I_{j,X}^{'}\supsetneq I_{j,X}$. By the maximality of $I_{j,X}$ for any
$j$ we have that $I_{j,X}^{'}\cap (X\setminus A_X)\neq \emptyset$, thus by lemma 3.2 iv) there exists $I\in S_{\phi}$ such that $I^{\ast}=X$ and
$I_{j,X}^{'}\cap I\neq\emptyset$. Since $I_{j,X}^{'}\cap A_X\neq\emptyset$, we conclude that $I_{j,X}^{'}\supsetneq I$, for any such
$I\in S_{\phi}$. We consider now a maximal disjoint subfamily of $(I_{j,X}^{'})_j$, denoted by $(I_{j_N,X}^{'})_N$, which still covers
$\cup_{j}I_{j,X}^{'}$. By the above construction we have that for every $N$, we can write $I_{j_N,X}^{'}=D_{j_N}\cup B_{j_N}$, where
$B_{j_N}=I_{j_N,X}^{'}\cap A_X$ and $D_{j_N}$ is a union of some of the elements $J$, of $S_{\phi}$ for which $J^{\ast}=X$. Obviously we have
$\cup_{N}B_{j_N}=A_X$ and each $B_{j_N}$ is a union of certain elements of the family $(I_{j,X})_j$. Now fix a $j_N$. For any $a\in(0,1)$ which will be chosen later, using
lemma 2.1, we construct a family $\ca_{\phi,j_N}^X$, of elements of $\ct$, all of which are contained in $B_{j_N}$, and such that
\begin{eqnarray}
\sum_{J\in \ca^X_{\phi,j_N}}\mi(J)=a\mi(B_{j_N}).  \label{eq3.13}
\end{eqnarray}

Define the function $g_{N,\phi,X}:B_{j_N}\ra\R^+$ by setting
\begin{eqnarray}
\begin{array}{lcl}
  g_{N,\phi,X}:=c^\phi_{N,X}, & \text{on} & \cup \ca^X_{\phi,j_N} \\
                \quad \; \quad \; \;\; :=0, &\text{on} & B_{j_N}\setminus\cup \ca^X_{\phi,j_N}
\end{array} \label{eq3.14}
\end{eqnarray}
where the constants $c^\phi_{N,X}$ and $\ga^\phi_{N,X}:=\mi(\cup \ca^X_{\phi,j_N})=a\mi(B_{j_N})$ satisfy
\begin{eqnarray}
\hspace{1cm}\left.\begin{array}{l}
        \int_{B_{j_N}}g_{N,\phi,X}d\mi=c^\phi_{N,X}\ga^\phi_{N,X}=\int\limits_{B_{j_N}}\phi d\mi \ \ \text{and} \\
        \int\limits_{B_{j_N}}g^p_{N,\phi,X} d\mi=(c^\phi_{N,X})^p\ga^\phi_{N,X}=\int\limits_{B_{j_N}}\phi^pd\mi,
      \end{array}\right\}, \;                                                                            \label{eq3.15}
\end{eqnarray}
It is easy to see that such choices for $c^\phi_{N,X}$ and $\ga^\phi_{N,X}$  are possible.

In fact (\ref{eq3.15}) give
\[
\ga^\phi_{N,X}=\left[\frac{\Big(\int\limits_{B_{j_N}}\phi d\mi\Big)^p}{\int\limits_{B_{j_N}}\phi^pd\mi}\right]^{1/(p-1)}\le\mi(B_{j_N}), \ \ \text{by H\"{o}lder's inequality}
\]
so we just need to define $\ga^\phi_{N,X}$, by the above equation, and choose $a$ so that
\[
a=\frac{\ga^\phi_{N,X}}{\mi(B_{j_N})}.
\]

At last we set $c^\phi_{N,X}=\dfrac{\int\limits_{B_{j_N}}\phi d\mi}{\ga^\phi_{N,X}}$.
Define now $g_{\phi}^{(1)}$ on $A_X=\cup _NB_{j_N}$ by $g_{\phi}^{(1)}(t)=g_{N,\phi,X}(t)$, for $t\in B_{j_N}$, for any $N$.
Note now that $g_{\phi}^{(1)}$may attain more than one positive values on $A_X$. It is easy then to see that there exists a common positive value,
denoted by $c^\phi_{X}$, and measurable sets $L_N \subseteq B_{j_N}$, such that if we define $g_{\phi}(t)=c^\phi_{X}$ for $t\in L_N$, and
$g_{\phi}(t)=0$, for $t\in B_{j_N}\setminus L_N$ and for any $N$, we still have that
$\int\limits_{B_{j_N}}g_\phi d\mi=\int\limits_{B_{j_N}}\phi d\mi=c^\phi_{X}\mi(L_N)$ and $\int\limits_{A_X}g^p_\phi d\mi=\int\limits_{A_X}\phi^p d\mi$.
For the construction of $L_N$ and $c^\phi_{X}$, we just need to find first the subsets $L_N$ of $B_{j_N}$ such that the first two of the integral equalities mentioned above is true, and this can be done for arbitrary $c^\phi_{X}$, since the space $(X,\mi)$ is nonatomic. Then we just need to find the constant $c^\phi_{X}$ for which the second integral equality is also true.
Note that for these choices of $L_N$ and $c^\phi_{X}$ we may not have
$\int\limits_{B_{j_N}}g^p_\phi d\mi=\int\limits_{B_{j_N}}\phi^p d\mi$, for every $N$, but the respective equality with $A_X$ in place of $B_{j_N}$ should be true.

Until now we have defined $g_\phi$ on $A_X$. We set now $g_\phi=\phi$ on $X\setminus A_X$. It is immediate then, by the construction of $g_\phi$,
that if $I\in \ct$ is such that $I\cap A_X\neq\emptyset$, and $I\cap (X\setminus A_X)\neq\emptyset$, we must have that
$\int\limits_{I}g_\phi d\mi=\int\limits_{I}\phi d\mi$. This is true since then $I$ can be written as a certain union of some subfamily of
$I_{j_N,X}^{'}$ and of some $J$, where $J$ is such that $J^{\ast}=X$. This last fact is true by the construction of the sets $I_{j_N,X}^{'}$, and because of lemma 3.3.

We continue then inductively and change the values of
$g_\phi$ on the sets $A_I$, for $I$, which is such that $I^{\ast}=X$, in the same way as was done before, but now working inside those $I$.
In this way we inductively define the function $g_\phi$ in all $X$, which obviously has the desired properties. Moreover the inequality
$\mi(\{\phi=0\}\cap A_I)\le \mi(\{g_{\phi}=0\}\cap A_I)$ is easily verified if we work as above in $B_{j_N}\cap \{\phi>0\}$ instead of
$B_{j_N}$. More precisely for the case of $I=X$ we define the family $\ca_{\phi,j_N}^X$, of elements of $\ct$, all of which are contained in $B_{j_N}$, by the relation
$\mi(\cup \ca^X_{\phi,j_N})=a\mi(B_{j_N}\cap\{\phi>0\})$, and define analogously $\ga^\phi_{N,X}$, now integrating on $B_{j_N}\cap\{\phi>0\}$.
Then we define in an analogous way $a$, that is we set $a=\frac{\ga^\phi_{N,X}}{\mi(B_{j_N}\cap\{\phi>0\})}$.
Now $\ga^\phi_{N,X}$ is less or equal than $\mi(B_{j_N}\cap\{\phi>0\})$, and by using this last fact we deduce that the zero set of $g_{\phi}$ in $A_X$,
increases in general, in relation to that of $\phi$ on the same set. The proof of our lemma is now completed.     \hs
\end{Proof}

Let now $(\phi_n)_n$ be an extremal sequence consisting of $\ct$-good functions and let $g_n=g_{\phi_n}$.\ We are now ready to prove the following
\begin{lem}\label{lem3.4}
With the above notation for an extremal $(\phi_n)_n$ sequence of $\ct$-good functions, we have that $\lim_n\mi(\{\phi_n=0\})=0$.
\end{lem}
\begin{Proof}
Fix $n\in\N$ and let $\phi=\phi_n$ and $g_\phi=g_{\phi_n}$ and $S=S_\phi$ the respective subtree of $\phi$.
We consider two cases:

i) $p\ge2$

We set $P_I=\dfrac{\int\limits_{A_I}\phi^pd\mi}{a_I}$, for every $I\in S_{\phi}$.
We obviously have $\ssum_{I\in S_\phi}a_IP_I=F$. We consider then the sum $\Sig_\phi=\ssum_{I\in S_\phi}\ga_IP_I$, where $\ga_I=\ga^\phi_I$ comes from lemma 3.3. More precisely, it should be true that $\int\limits_{A_I}\phi d\mi=\ga_Ic_I$ and $\int\limits_{A_I}\phi^p d\mi=\ga_Ic^p_I$, for a suitable constant
$c_I=c^{\phi}_I$. Obviously $0\leq\ga_I\leq a_I=\mi(A_I)$, so we must have that
\begin{align*}
\Sig_\phi&=\sum_{I\in S_\phi}\ga_I\frac{\int\limits_{A_I}\phi^p}{a_I}=
\sum_{I\in S_\phi}\ga_I\frac{\ga_I\cdot c^p_I}{a_I}
=\sum_{I\in S_\phi}\ga^2_I\frac{c^p_I}{a_I}=\sum_{I\in S_\phi}
\frac{\ga^2_Ia^{p-2}_Ic^p_I}{a^{p-1}_I} \\
&\overset{p\ge2}{\ge}\sum_{I\in S_\phi}\frac{(\ga_Ic_I)^p}{a^{p-1}_I}=
\sum_{I\in S_\phi}\frac{\Big(\int\limits_{A_I}\phi\Big)^p}{a_I^{p-1}}.
\end{align*}
From the first inequality in (4.20) in \cite{Mel1}, and since $\phi_n$ is extremal we have that the last sum in the above inequality tends to $F$, as $\phi$ moves along $(\phi_n)_n$. We conclude that
\begin{eqnarray}
\sum_{I\in S_{\phi}}\ga_IP_I\approx F,   \label{eq3.16}
\end{eqnarray}
since  $\Sig_\phi\le F$.
Consider now for every $R>0$ and every $\phi$ the following set
\[
S_{\phi,R}=\cup\{A_I=A(\phi,I):\;I\in S_\phi, \ \ P_I<R\}.
\]
For every $I\in S_\phi$ such that $P_I<R$ we have that $\int\limits_{A_I}\phi^p<Ra_{I}$. Summing for all such $I$ we obtain
\begin{eqnarray}
\int_{S_{\phi,R}}\phi^pd\mi<R\mi(S_{\phi,R}).  \label{eq3.17}
\end{eqnarray}
Additionally we have that
\begin{eqnarray}
\left|\sum_{I\in S_\phi\atop P_I\ge R}a_IP_I-F\right|=\int_{S_{\phi,R}}\phi^p d \mi, \ \ \text{and}  \label{eq3.18}
\end{eqnarray}
\begin{eqnarray}
\sum_{I\in S_\phi\atop P_I<R}\ga_IP_I\le\sum_{I\in S_\phi\atop P_I<R}a_IP_I=\int_{S_{\phi,R}}\phi^p d\mi.  \label{eq3.19}
\end{eqnarray}
From (\ref{eq3.16}) and (\ref{eq3.19}) we have that
\begin{eqnarray}
\underset{\phi}{\lim\sup}\left|\sum_{I\in S_\phi\atop P_I\ge R}\ga_IP_I-F\right|\le\lim_\phi\int_{S_{\phi,R}}\phi^pd\mi,  \label{eq3.20}
\end{eqnarray}
where we have supposed that the last limit exists (in the opposite case we just pass to a subsequence of $(\phi_n)_n$).\ From (\ref{eq3.18}) and (\ref{eq3.20}) we conclude that
\begin{eqnarray}
\underset{\phi}{\lim\sup}\sum_{I\in S_\phi\atop P_I\ge R}(a_I-\ga_I)P_I\le 2\lim_\phi\int_{S_{\phi,R}}\phi^pd\mi.  \label{eq3.21}
\end{eqnarray}
By using now theorem \ref{thm3.3} we have that
\[
\lim_\phi\int_{K_\phi}(\cm_\ct\phi)^pd\mi=\oo_p(f^p/F)^p\lim_\phi\int_{K\phi}\phi^pd\mi,
\]
whenever the limits exist, and $K_\phi$ be a union of pairwise disjoint elements of $S_\phi$ (the conditions of theorem \ref{thm3.3} are satisfied because of the boundedness of the sequences mentioned there).

Now for a fixed $R>0$, $S_{\phi,R}$ is a union of sets of the form $A_I$, for certain $I\in S_\phi$. Each $A_I$ can be written in view of lemma \ref{lem3.2} as $A_I=I\setminus\bigcup\limits_{J\in S_\phi, J^{\ast}=I}J$.\ Using then a diagonal argument and passing if necessary to a subsequence we can suppose that
\begin{eqnarray}
\lim_\phi\int_{S_{\phi,R}}(\cm_\ct\phi)^pd\mi=\oo_p(f^p/F)^p\lim_\phi\int_{S_{\phi,R}}\phi^p,
\label{eq3.22}
\end{eqnarray}
by applying theorem 3.3 as mentioned above.
Since now $\cm_\ct\phi(t)\ge f$, for every $t\in X$, we have that
\begin{eqnarray}
\lim_\phi\int_{S_{\phi,R}}(\cm_\ct\phi)^pd\mi\ge(\underset{\phi}{\lim\sup}\mi(S_{\phi,R}))f^p, \label{eq3.23}
\end{eqnarray}
and because of (\ref{eq3.17}) we have that
\begin{eqnarray}
\lim_\phi\int_{S_{\phi,R}}\phi^pd\mi\le\underset{\phi}{\lim\sup} R\mi(S_{\phi,R}),  \label{eq3.24}
\end{eqnarray}
for any $R>0$.\ Combining the last two relations (in view of (\ref{eq3.22})) we obtain that
\begin{eqnarray}
f^p(\underset{\phi}{\lim\sup}\mi(S_{\phi,R}))\le R\oo_p(f^p/F)^p\cdot
(\underset{\phi}{\lim\sup}\mi(S_{\phi,R})),  \label{eq3.25}
\end{eqnarray}
so by choosing $R>0$ suitable small depending only on $f,F$ we have that
\begin{eqnarray}
\underset{\phi}{\lim\sup}\mi(S_{\phi,R})=0.  \label{eq3.26}
\end{eqnarray}
Using now (\ref{eq3.21}) and (\ref{eq3.24}) we obtain, for this $R$, that
\[
R\underset{\phi}{\lim\sup}\sum_{I\in S_\phi\atop P_I\ge R}(a_I-\ga_I)\le2\lim_\phi\int_{S_{\phi,R}}\phi^pd\mi\le2R\lim_\phi\mi(S_{\phi,R})=0
\]
Thus
\begin{eqnarray}
\lim_\phi\sum_{I\in S_\phi\atop P_I\ge R}(a_I-\ga_I)=0. \label{eq3.27}
\end{eqnarray}
Since now $\ssum_{I\in S_\phi}a_I=1$, $\mi(S_{\phi,R})=\ssum_{I\in S_\phi\atop P_I<R}a_I$ we easily obtain from (\ref{eq3.27}) that:
\begin{align*}
&\lim_\phi\left[1-\mi(S_{\phi,R})-\sum_{I\in S_\phi\atop P_I\ge R}\ga_I\right]=0 \;\Rightarrow \\
&\lim_\phi\sum_{I\in S_\phi\atop P_I\ge R}\ga_I=1 \;\Rightarrow \\
&\lim_\phi\sum_{I\in S_\phi}\ga_I=1 \;\Rightarrow \\
&\lim_\phi\sum_{I\in S_\phi}(a_I-\ga_I)=0.
\end{align*}
Thus we must have that
\[
\mi(\{\phi=0\})\le\mi(\{g_\phi=0\})=\sum_{I\in S_\phi}(a_I-\ga_I)\overset{\phi}{\longrightarrow}0.
\]
Lemma 3.2 is proved in the first case.

ii) The case $1<p<2$ is treated in a similar way:

Here we define $P_I=\dfrac{\int\limits_{A_I}\phi^p}{a^{p-1}_I}$ and prove in the same manner that
\[
\lim_\phi\sum_{I\in S_\phi}(a^{p-1}_I-\ga^{p-1}_I)P_I=0.
\]
Using then the inequality $x^{q}-y^{q}>q(x-y)$, which holds for $1>x>y$ and $0<q<1$, we conclude that:
\begin{align*}
&\lim_\phi\sum_{I\in S_\phi}(a_I-\ga_I)=0  \;\Rightarrow \\
&\lim_\phi\mi(\{g_\phi=0\})=0 \;\Rightarrow \\
&\lim_\phi\mi(\{\phi=0\})=0,
\end{align*}
and by this we end the proof of lemma 3.4.  \hs
\end{Proof}

Suppose now that $(\phi_n)_n$ is extremal. For every $\phi\in\{\phi_n,\;n=1,2,\ld\}$
we define $g'_\phi:x\ra\R^+$ by  $g'_\phi(t)=c^\phi_I$, $t\in A_I$ for $I\in S_\phi$, that is we ignore the zero values of $g_\phi$. Then we easily see because of lemma 3.4 that
\[
\lim_\phi\int_Xg'_\phi d\mi=f, \ \ \lim_\phi\int_X(g'_\phi)^pd\mi=F \ \ \text{and}
\]
\begin{eqnarray}
\lim_\phi\int_X|g_\phi-g'_\phi|^pd\mi=0.  \label{eq3.28}
\end{eqnarray}
Obviously also by lemma 3.3 we have that
\begin{eqnarray}
Av_{I}(g_\phi)=Av_{I}(\phi),  \label{eq3.29}
\end{eqnarray}
for every $I\in S_\phi$.
From (\ref{eq3.29}) we have that $\cm_\ct g_\phi\ge\cm_\ct\phi$ on $X$ $\Rightarrow$ $\dis\lim_\phi\int\limits_X(\cm_\ct g_\phi)^pd\mi=F\oo_p(f^p/F)^p$, in view of (\ref{eq3.15}) and theorem \ref{thm2.1}.

Since $\int\limits_Xg_\phi d\mi=f$, $\int\limits_X(g_\phi)^pd\mi=F$ we have that $(g_\phi)_\phi$ is an extremal sequence. Suppose now that we have proved the following two equalities
\begin{eqnarray}
\lim_\phi\int_X|g'_\phi-\phi|^pd\mi=0,  \label{eq3.30}
\end{eqnarray}
and
\begin{eqnarray}
\lim_\phi\int_X|\cm_\ct g_\phi-cg_\phi|^pd\mi=0,\ \  \text{for} \ \  c=\oo_p(f^p/F). \label{eq3.31}
\end{eqnarray}
Then because of (\ref{eq3.28}) we would have that
\begin{align*}
&\lim_\phi\int_X|\phi-g_\phi|^pd\mi=0\overset{(3.31)}{\Rightarrow} \\
&\lim_\phi\int_X|\cm_\ct\phi-c\phi|^pd\mi=0
\end{align*}
which is the result we need to prove.
We proceed to the proof of (\ref{eq3.30}) and (\ref{eq3.31}).\vspace*{0.2cm} \\
\begin{lem}\label{lem3.5}
With the above notation
\[
\lim_\phi\int_X|\cm_\ct g_\phi-cg_\phi|^qd\mi=0.
\]
\end{lem}
\begin{Proof}
We recall that $c=\oo_p(f^p/F)$.
We set for each $\phi\in\{\phi_n$, $n=1,2,\ld\}$
\[
\De_\phi=\{t\in X:\cm_\ct g_\phi(t)> cg_\phi(t)\}
\]
It is obvious by passing if necessary to a subsequence that
\begin{eqnarray}
\lim_\phi\int_{\De_\phi}(\cm_\ct g_\phi)^pd\mi\ge\oo_p(f^p/F)^p\lim_\phi\int_{\De_\phi}g^p_\phi d\mi.
\label{eq3.32}
\end{eqnarray}
We consider now for every $I\in S_\phi$ the set $(X\setminus\De_\phi)\cap A_I$. We distinguish two cases:

(i) $Av_I(\phi)=y_I>cc_I^{\phi}$, where $c_I^{\phi}$ is the positive value of $g_\phi$ on $A_I$ (if it exists). Then because of (\ref{eq3.29}) we have that
$\cm_{\ct}g_{\phi}(t)\geq Av_I(g_{\phi})=Av_I(\phi)>cc_I^{\phi}\geq cg_{\phi}(t)$, for each $t\in A_I$. Thus $(X\setminus\De_\phi)\cap A_I=\emptyset$
in this case. We study now the second case.

(ii) $y_I\leq cc_I^{\phi}$. Let now $t\in A_I$ with $g_{\phi}(t)>0$, that is $g_{\phi}(t)=c_I^{\phi}$. We prove that for each such $t$ we have
$\cm_{\ct}g_{\phi}(t)\leq cg_{\phi}(t)=cc_I^{\phi}$. Suppose now that for some $t$ we have the opposite inequality. Then there exists $J_t$ such that
$t\in J_t$ and $Av_{J_t}(g_{\phi})>cc_I^{\phi}$. Then one of the following subcases holds

(a) $J_t\subseteq A_I$. Then by the form of $g_{\phi}/A_I$ (equals $0$ or $c_I^{\phi}$), we have that $Av_{J_t}(g_{\phi})\leq c_I^{\phi}<cc_I^{\phi}$,
which is a contradiction, since $c>1$. Thus this case is excluded.

(b) $J_t$ is not a subset of $A_I$. Then in this subcase two more subcases can occur.

$b_1)$ $J_t\subseteq I$ and $J_t$ contains properly an element of $S_{\phi}$, $J'$, for which $(J')^{\ast}=I$. Since now (ii) holds, $t\in J_t$ and
$Av_{J_t}(g_{\phi})>cc_I^{\phi}$, we must have that $J'\subsetneq J_t \subsetneq I$. We choose now an element of $\ct$, $J'_t\subsetneq I$, which contains
$J_t$, with maximum value on the average $Av_{J'_t}(\phi)$. Then by it's choice we have that for each $K\in \ct$ such that $J'_t\subseteq K \subsetneq I$
the following holds $Av_{K}(\phi)\leq Av_{J'_t}(\phi)$. Since now $I\in S_{\phi}$ and $y_I=Av_I(\phi)\leq cc_I^{\phi}$ by lemma 3.1 and the choice of $J'_t$ we have that $Av_K(\phi)< Av_{J'_t}(\phi)$ for every $K\in \ct$ such that $J'_t\subsetneq K$. So again by lemma 3.1 we conclude that $J'_t\in S_{\phi}$. But this is impossible since $J'\subsetneq J'_t\subsetneq $, $J',I\in S_{\phi}$ and $(J')^{\ast}=I$. We turn now to the last subcase.

$b_2)$ $I\subsetneq J_t$. Then by an application of lemma 3.3 we have that $Av_{J_t}(\phi)=Av_{J_t}(g_{\phi})>cc_I^{\phi}\geq y_I=Av_I(\phi)$ which is
impossible by lemma 3.1, since $I\in S_{\phi}$.

In any of the two cases $b_1)$ and $b_2)$ we have proved that we have $(X\setminus\De_\phi)\cap A_I= A_I\setminus (g_{\phi}=0)$, while we showed that in case
(i), $(X\setminus\De_\phi)\cap A_I=\emptyset$.

Since $\bigcup_{I\in S_{\phi}}A_I\approx X$ we conclude by lemma 3.4 and the above discussion that
$X\setminus\De_\phi\approx \big(\bigcup_{I\in S_{1,\phi}}A_I\big)\setminus E_{\phi}$, where $\mi(E_{\phi})\rightarrow 0$ and $S_{1,\phi}$ is a subset of the
subtree $S_{\phi}$. Since now each $A_I, I\in S_{1,\phi}\subseteq S_{\phi} $ is written, by lemma 3.2, as a difference set of unions of elements of $S_{\phi}$,
and theorem 3.3 holds for such unions, we conclude by a diagonal argument and by passing if necessary to a subsequence, that
\[
\lim_\phi\int_{\cup A_I\atop I\in S_{1,\phi}}(\cm_\ct\phi)^pd\mi=\oo_p(f^p/F)^p\cdot\lim_\phi\int_{\cup A_I\atop I\in S_{1,\phi}}\phi^pd\mi, \ \ \text{and since}
\]
\[
\mi(E_\phi)\ra0\Longrightarrow\lim_\phi\int_{X\setminus\De_\phi}(\cm_\ct\phi)^pd\mi=\oo_p(f^p/F)^p\lim_\phi\int_{X\setminus\De_\phi}
\phi^pd\mi.
\]
Because now of the relation $\cm_\ct g_{\phi}\geq \cm_\ct\phi$ ,which holds $\mi$-almost everywhere on $X$, we have as a result that
\begin{eqnarray}
\lim_\phi\int_{X\setminus\De_\phi}(\cm_\ct g_\phi)^pd\mi\ge\oo_p(f^p/F)^p\lim_\phi\int_{X\setminus\De_\phi}g^p_\phi d\mi.
\label{eq3.33}
\end{eqnarray}
Adding the relations (3.32) and (3.33) we have obtained that $\dis\lim_\phi\int\limits_X(\cm_\ct g_{\phi})^pd\mi\ge\oo_p(f^p/F)^pF$, which in fact is an equality
since $(g_{\phi})$ is an extremal sequence. So we must have equality in both (3.32) and (3.33). By using then the elementary inequality $x^p-y^p>(x-y)^p$ which holds for every $x>y>0$ and $p>1$, in view of the inequality $\cm_\ct g_\phi>cg_\phi$, which holds on $\De_\phi$, we must have that
\begin{eqnarray}
\lim_\phi\int_{\De_\phi}|\cm_\ct g_\phi-cg_\phi|^pd\mi=0.
\label{eq3.34}
\end{eqnarray}
Similarly for $X\setminus \De_\phi$. That is
\begin{eqnarray}
\lim_\phi\int_{\De_\phi}|\cm_\ct g_\phi-cg_\phi|^pd\mi=0.
\label{eq3.35}
\end{eqnarray}
Adding (3.34) and (3.35) we derive $\lim_\phi||\cm_\ct g_\phi-cg_\phi||_{L^p}=0$, and by this we end the proof of our lemma. \hs
\end{Proof}

We now proceed to
\begin{lem}\label{lem3.4}
Under the above notation (\ref{eq3.30}) is true.
\end{lem}
\begin{Proof}
We just need to prove that
\begin{eqnarray}
\lim_\phi\int_{\{g'_\phi\le\phi\}}[\phi^p-(g'_\phi)^p]d\mi=0.  \label{eq3.36}
\end{eqnarray}
Then since
\[
\lim_\phi\int_{\{g'_\phi\le\phi\}}[\phi^p-(g'_\phi)^p]d\mi=\lim_\phi\int_{\phi\le g'_\phi}
[(g'_\phi)^p-\phi^p], \ \ \text{and} \ \ p>1
\]
we have the desired result, in view of the inequality $(x-y)^p<x^p-y^p$, which holds for $0<y<x$ and $p>1$.

We use the inequality
\begin{eqnarray}
\hspace*{0.5cm}t\le\frac{t^p}{p}+\frac{1}{q}, \;\text{for every}\; t>0 \; \text{where}\; p,q>1\;\text{such that}\; \frac{1}{p}+\frac{1}{q}=1,  \label{eq3.37}
\end{eqnarray}
For any $I\in S_{\phi}$ we set
\begin{align*}
&\De^{(1)}_{I,\phi}=\{g'_\phi\le\phi\}\cap A(\phi,I) \\
&\De^{(2)}_{I,\phi}=\{\phi< g'_\phi\}\cap A(\phi,I).
\end{align*}
Because of (\ref{eq3.37}), if we write $c_{I,\phi}$ instead of $c^{\phi}_I$ and suppose that $c_{I,\phi}>0$, we have that
\[
\frac{1}{c_{I,\phi}}\phi(x)\le\frac{1}{p}\frac{1}{c^p_{I,\phi}}\phi^p(x)+\frac{1}{q}, \ \ \text{for every} \ \ x\in A_I=A(\phi,I).
\]
Integrating over $\De^{(1)}_{I,\phi}$, and $\De^{(2)}_{I,\phi}$ we have that
\[
\frac{1}{c_{I,\phi}}\int_{\De^{(j)}_{I,\phi}}\phi d\mi\le\frac{1}{p}\frac{1}{c^p_{I,\phi}}\int_{\De^{(j)}_{I,\phi}}\phi^pd\mi+\frac{1}{q}
\mi(\De^{(j)}_{I,\phi}), \ \ \text{for} \ \ j=1,2, \ \ I\in S_{\phi}
\]
which gives
\[
c^{p-1}_{I,\phi}\int_{\De^{(j)}_{I,p}}\phi d\mi\le\frac{1}{p}\int_{\De^{(j)}_{I,\phi}}\phi^pd\mi+\frac{1}{q}\mi(\De^{(j)}_{I,\phi})c^p_{I,\phi}.
\]
Note that the last inequality is satisfied even if $c_{I,\phi}=0$. Summing the above for $I\in S_\phi$ we obtain
\begin{eqnarray}
\sum_{I\in S_\phi}c^{p-1}_{I,\phi}\int_{\De^{(j)}_{I,\phi}}\phi d\mi\le\frac{1}{p}\int_{\bigcup\limits_I\De^{(j)}_{I,\phi}}\phi^p d\mi+
\frac{1}{q}\sum_{I\in S_\phi}\mi(\De^{(j)}_{I,\phi})c^p_{I,\phi}, \label{eq3.38}
\end{eqnarray}
for $j=1,2$, thus by adding the above two inequalities we conclude that
\begin{eqnarray}
\sum_{I\in S_\phi}c^{p-1}_{I,\phi}\int_{A(\phi,I)}\phi d\mi\le\frac{1}{p}F+\frac{1}{q}\sum_{I\in S_\phi}\mi(A(\phi,I))c^p_{I,\phi}.  \label{eq3.39}
\end{eqnarray}
The left hand side of (\ref{eq3.39}) is equal to
\[
\sum_{I\in S_\phi} c^{p-1}_{I,\phi}(c_{I,\phi}\ga^\phi_I)=\sum_{I\in S_\phi}\ga^\phi_Ic^p_{I,\phi}=\int_Xg^p_\phi d\mi
\]
while the right hand side is equal to $\dfrac{1}{p}F+\dfrac{1}{q}\int\limits_X(g'_\phi)^pd\mi$. In the limit we have equality on (\ref{eq3.39}), because of (\ref{eq3.28}). This gives equality on (\ref{eq3.38}) for $j=1,2$ in the limit. Thus for $j=1$ we have that
\[
\sum_{I\in S_\phi}c^{p-1}_{I,\phi}\int_{\De^{(1)}_{I,\phi}}\phi d\mi\approx\frac{1}{p}\sum_{I\in S_\phi}\int_{\De^{(1)}_{I,\phi}}\phi^pd\mi+\frac{1}{q}\sum_{I\in S_\phi} c^p_{I,\phi}\mi(\De^{(1)}_{I,\phi})\Rightarrow
\]
\begin{eqnarray}
\int_{\{g'_\phi\le\phi\}}\phi(g'_\phi)^{p-1}d\mi\approx\frac{1}{p}\int_{\{g'_\phi\le\phi\}}
\phi^pd\mi+\frac{1}{q}\int_{\{g'_\phi\le\phi\}}(g'_\phi)^pd\mi.  \label{eq3.40}
\end{eqnarray}
We set
\[
t_\phi=\left(\int_{\{g'_\phi\le\phi\}}\phi^p d\mi\right)^{1/p}, \ \ S_\phi=\left(\int_{\{g'_\phi\le\phi\}}(g'_\phi)^pd\mi\right)^{1/p}.
\]
Then
\[
\int_{\{g'_\phi\le\phi\}}\phi(g'_\phi)^{p-1}d\mi\le t_\phi\cdot S^{p-1}_\phi, \ \ \text{so (\ref{eq3.40}) gives:}
\]
\[
\frac{1}{p}t^p_\phi+\frac{1}{q}S^p_\phi\underset{\phi}{\le}t_\phi\cdot S^{p-1}_\phi,
\]
so as a result we have because of (\ref{eq3.37}) that
\[
\frac{1}{p}t^p_\phi+\frac{1}{q}S^p_\phi\underset{\phi}{\approx}t_\phi\cdot S^{p-1}_\phi.
\]
Since now in (\ref{eq3.37})  we have equality only for $t=1$, and $t_\phi,S_\phi$ are bounded we conclude that
\[
\frac{t^p_\phi}{S^p_\phi}\underset{\phi}{\longrightarrow}1\Rightarrow t^p_\phi-S^p_\phi\overset{\phi}{\longrightarrow}0\Rightarrow\int_{\{g'_\phi\le\phi\}}
[\phi^p-(g'_\phi)^p]d\mi\overset{\phi}{\longrightarrow}0,
\]
which is (\ref{eq3.36}). The proof of our lemma is now completed.   \hs
\end{Proof}
We have thus completed the proof of theorem A. We should also mention that since $\ct$-good functions include $\ct$-step functions, in the case of $\R^n$, where the Bellman function is given by (\ref{eq1.4}) for a fixed dyadic cube $Q$, we obtain the result in theorem A for every sequence of Lesbesgue measurable functions $(\phi_n)_n$.\ In general in all interesting cases we do not need the hypothesis for $\phi_n$ to be $\ct$-good since $\ct$-simple functions are dense on $L^p(X,\mi)$.

Department of Mathematics,
National and Kapodistrian University of Athens
Panepistimioupolis, GR 157 84, Athens,
Greece

\end{document}